\documentclass[11pt]{amsart}
\usepackage{preamble}
\usepackage{placeins}

\begin{document}

\title[Full-window branch discovery and EnKF continuation]{
Full-window branch discovery and loss-selected EnKF continuation for data assimilation
}

\begin{abstract}
We develop a framework for offline full-window branch discovery, optionally followed by online continuation with an ensemble Kalman filter (EnKF).
Three mechanisms drive the branch search: adjoint path-kernel (APK) differentiation balances kernel differentiation and correction-stabilized path perturbation, shifting the optimization from exploration to exploitation; an optimized Gaussian initial law broadens the search over initial-state basins; and loss-weighted mixing across independent runs recombines successful path components.
We may then select an interior state using a local loss and continue online with an EnKF.
In \(40\)-dimensional Lorenz--96 experiments, the mean offline path RMSE of APK is \(4.3\) times smaller than that of population weak-\(\mathrm{4D\text{-}Var}_x\).
The resulting APK--EnKF method has a mean online RMSE \(64\) times smaller than that of ordinary EnKF.
\end{abstract}

\maketitle

\section{Introduction}
\label{s:intro}

\subsection{Main results}
\label{s:mainResults}

Long-window data assimilation for chaotic systems presents two distinct difficulties.
First, ordinary path perturbations can become exponentially unstable, making useful derivatives difficult to compute.
Second, sparse and partial observations can admit several dynamically plausible paths, so independent optimizations may settle into branch-dependent local minima even when their gradients are stable.

Our offline branch-discovery framework combines three mechanisms.
The first is APK differentiation, which combines kernel differentiation with correction-stabilized path perturbation \cite{dud,apk}.
Early in the optimization, prescribed large values of \(\alpha\) and the diffusion scales support exploration through random paths and kernel differentiation.
As these prescribed scales decrease, the optimization shifts toward pathwise local exploitation.
The second mechanism represents the physical initial state by an optimized Gaussian law rather than a deterministic control, broadening the search over initial-state basins.
The third mechanism periodically mixes independently optimized path reconstructions across physical time and state coordinates.
Loss-dependent weights computed from the offline assimilation window recombine successful path components into a proposal that replaces the current worst population member.
This allows information discovered by different runs to be combined when no single run has yet followed a good branch throughout the window.

After the offline branch search identifies a reconstructed path, an optional continuation step uses a local loss to choose an interior restart state without reference to the truth.
The EnKF remains a separate online component: it is initialized from this state, replayed through the remainder of the offline assimilation window, and then continued online, provided that the selected branch is locally trackable.

We test these ideas in twenty independently seeded realizations of a \(40\)-dimensional Lorenz--96 system, observing one coordinate and hiding the next four, with independent observation noise of standard deviation \(0.3\).
APK has a mean offline path RMSE of \(0.698\), compared with \(3.02\) for population weak-\(\mathrm{4D\text{-}Var}_x\), and a restart-state RMSE of \(0.180\), compared with \(2.56\).
The APK-initialized EnKF remains finite in all experiments; in the thirteen experiments in which ordinary EnKF also remains finite, their mean online RMSEs are \(0.0959\) and \(6.12\), respectively.
A separate nonlinear-observation experiment with squared observations gives the same offline ordering and produces an APK--EnKF online RMSE below \(0.5\) in fifteen of twenty experiments.

\subsection{Literature review}
\label{s:review}

Data assimilation (DA) combines imperfect dynamical models with noisy, partial observations \cite{DA_Kalnay,DA_AST}.
One can distinguish sequential filtering methods from four-dimensional variational assimilation (4D-Var).
Filtering methods propagate a state estimate and update it as observations arrive \cite{filter_ABN,filter_NMX,filter_LDZ,filter_Yau,Sanz-Alonso2025}.
Once an EnKF is initialized near the correct branch, it can provide efficient local tracking.
However, its causal Gaussian updates are not designed to compare complete branch histories retrospectively, and a discarded branch cannot be recovered using later observations without an additional smoothing or full-window procedure.
Accordingly, experiments using accurate or oracle initialization supply the branch in advance and mainly test local tracking.
Our framework instead discovers a branch from the offline assimilation window and then optionally initializes an EnKF for online continuation.

By contrast, 4D-Var uses the entire assimilation window to optimize a trajectory and, possibly, dynamical parameters.
Strong-constraint 4D-Var takes the initial condition as its control and enforces the model exactly.
In the forcing formulation of weak-constraint 4D-Var, denoted \(\mathrm{4D\text{-}Var}_q\), the controls also include model-error increments \(q_n\), and the states are generated sequentially by
\[
x_{n+1}=\mathcal M_n(x_n)+q_n.
\]
In both formulations, perturbations of the initial condition or early model-error increments propagate through the remaining window via long products of tangent maps and may grow exponentially in chaotic systems \cite{4dvar_DT,4dvar_TC}.

Weak-constraint 4D-Var also has a state formulation, denoted \(\mathrm{4D\text{-}Var}_x\), in which \(x_0,\ldots,x_N\) are optimized directly and the local model defects \(x_{n+1}-\mathcal M_n(x_n)\) are penalized \cite{Tremolet2006,XuAnitescu2016}.
The gradient with respect to each state then involves only adjacent time intervals and avoids tangent products across the full window.
Nevertheless, the deterministic objective can contain competing local minima, while the local path updates may transfer information slowly across the window \cite{Cox2015,Brocker2017}.

In this paper, the complete branch-discovery framework combines three mechanisms: APK differentiation provides stable full-window gradients and kernel-based exploration, an optimized Gaussian initial law broadens the search over initial-state basins, and population mixing recombines successful portions of different full-window paths.
The resulting branch search remains an offline smoothing procedure that uses the observations \(y_{0:N}\) over the assimilation window to estimate the state path \(x_{0:N}\).

The gradient in our variational solver is computed by the adjoint path-kernel method, which combines path perturbation and kernel differentiation.
Path perturbation fixes the Brownian path and differentiates the resulting state path, as in stochastic-gradient, backpropagation, and adjoint methods \cite{eyink2004ruelle,lucarini_linear_response_climate,Caflisch2021}; for chaotic dynamics, however, the perturbation and estimator variance can grow exponentially with time.
Kernel differentiation instead fixes the realized state path and differentiates the probability kernel of the noise \cite{CM44,MalliavinBook,Rubinstein1989,Reiman1989,Glynn1990}.
Together with the artificial process noise, it supports exploration through random paths and avoids propagating tangent Jacobians, but its variance increases as the noise decreases.
The path-kernel method combines these mechanisms by damping unstable path perturbations and recovering the removed contribution through kernel differentiation \cite{dud}.
Its adjoint form, APK, supplies the stable full-window gradients used here \cite{apk}.

The linear correction and population mixing are related to previous methods but serve different roles here.
Feedback corrections have long been used in nudging and continuous data assimilation to relax model trajectories toward observations \cite{AurouxBlum2008,AzouaniOlsonTiti2014}.
Here, the feedback target is instead an optimized full-state path \(x^*\), and the correction stabilizes path perturbations within a penalized full-window optimization.
Genetic searches have also used populations and random crossover to optimize weak-constraint trajectories \cite{Ahrens1999}.
Our mixing instead assigns credit at each physical time and state coordinate and recombines locally successful path portions.

After solving the full-window problem and finding the correct branch, we select where within the offline assimilation window the EnKF continuation should begin.
Existing hybrid 4D-Var/EnKF methods use prescribed cycle times, sometimes at the middle of the window, while analysis-time studies compare predefined beginning, middle, and end times \cite{ZhangZhang2012E4DVar,PoterjoyZhang2015,ZhaoWangLiu2011}.
Ensemble smoothers likewise use a prescribed lag or window \cite{BocquetSakov2014IEnKS}.
Here, the restart time is instead selected separately for each reconstructed path by minimizing a truth-free local loss over admissible interior times.
The EnKF is initialized from the corresponding state, replayed through the remaining observations in the offline assimilation window, and then continued online.
In our experiments, this typically selects a more accurate restart state with lower RMSE.

\subsection{Structure of the paper}
\label{s:structure}

\Cref{s:notations} introduces the discrete filtering and path-space notation.
\Cref{s:theory} develops the complete branch-discovery framework, including APK gradients with an optimized Gaussian initial law, the prescribed transition from kernel-based exploration to pathwise exploitation, offline path reconstruction, population mixing, and loss-selected optional EnKF continuation.
\Cref{s:numeric} evaluates offline reconstruction against population weak-\(\mathrm{4D\text{-}Var}_x\) and EnKF continuation against ordinary EnKF in \(40\)-dimensional Lorenz--96 experiments with linear and squared observations.
\Cref{s:conclusion} summarizes the conclusions and future work on adaptive exploration in time and space.

\section{Notations and Preparations}
\label{s:notations}

\label{s:filter_notation}

We write \(x_n\in\R^M\) for the model state and \(y_n\in\R^m\) for the observation at time \(n\), with \(m\le M\).
For a path, we use the shorthand
\[
x_{0:N}:=(x_0,\ldots,x_N),
\qquad
y_{0:N}:=(y_0,\ldots,y_N).
\]
The observation map is denoted by \(\phi:\R^M\to\R^m\).
We fix a time step \(\Dt\), write \(T=N\Dt\), and call \([0,T]\) the offline assimilation window.
We use the It\^o/Euler convention
\[
\DB_n:=B_{(n+1)\Dt}-B_{n\Dt},
\qquad
\DB_n\iid \cN(0,\Dt I).
\]
The natural filtration of the sampled model path is
\[
\cF_n:=\sigma(x_0,\DB_0,\ldots,\DB_{n-1}).
\]
Thus a causal filter estimates \(x_n\) using only observations up to time \(n\).
By contrast, APK uses all observations \(y_{0:N}\) in the offline assimilation window to estimate the whole state path \(x_{0:N}\).
This is the smoothing or path-space viewpoint, and it is the reason later observations can correct earlier state estimates.

Throughout this paper, the online filter is applied only after this offline initialization.
The expensive part is the path-space optimization over the offline assimilation window; the online component is a fast local continuation.

We use the convention that a subscript \(n\) on a function or derivative means evaluation along the sampled path; for example,
\[
F_n=F(x_n),
\qquad
\nabla F_n=\nabla F(x_n),
\qquad
\phi_n=\phi(x_n).
\]
We assume throughout that the displayed derivatives and expectations are well-defined.

\section{Full-window branch-discovery framework}
\label{s:theory}

This section develops the complete branch-discovery framework.
We first derive APK gradients with an optimized random initial condition and explain the prescribed transition from kernel-based exploration to pathwise exploitation.
We then present offline path reconstruction, population mixing, and optional loss-selected EnKF continuation.

\subsection{APK with an optimized random initial condition}
\label{s:pk_review}

This subsection states the finite-time path-kernel formulas used by APK and extends them to an optimized Gaussian initial law.
The finite-time identity is proved in \cite{dud}, and its adjoint form is derived in \cite{apk}; we state only the formulas needed here.

\subsubsection{Finite-time path-kernel identity}
For simplicity, we use scalar diffusion.
Let \(\boldsymbol{\gamma}\in\R^{N_\gamma}\) collect all optimization parameters.
Let \(\boldsymbol{\gamma}=0\) denote the parameter value at which the gradient is evaluated, and write
\[
\delta A^{\boldsymbol{\gamma}}
:=\left.\frac{\partial A^{\boldsymbol{\gamma}}}{\partial\boldsymbol{\gamma}}\right|_{\boldsymbol{\gamma}=0}.
\]
For \(A^{\boldsymbol{\gamma}}\in\R^q\), \(\delta A^{\boldsymbol{\gamma}}\in\R^{q\times N_\gamma}\) is its parameter Jacobian.
We reserve \(\nabla\) for spatial derivatives and omit the superscript \(\boldsymbol{\gamma}=0\) from baseline quantities.

Consider the time-discretized SDE
\[
x_{n+1}^{\boldsymbol{\gamma}}
=x_n^{\boldsymbol{\gamma}}+F^{\boldsymbol{\gamma}}(x_n^{\boldsymbol{\gamma}})\Dt
+\sigma^{\boldsymbol{\gamma}}(x_n^{\boldsymbol{\gamma}})\DB_n,
\qquad \DB_n\iid\cN(0,\Dt I).
\]
In \(\delta F_n^{\boldsymbol{\gamma}}:=\delta F^{\boldsymbol{\gamma}}(x_n)\) and
\(\delta\sigma_n^{\boldsymbol{\gamma}}:=\delta\sigma^{\boldsymbol{\gamma}}(x_n)\), the spatial input \(x_n\) is held fixed.
Let \(\alpha_n\) be an \(\cF_n\)-adapted scalar schedule independent of \(\boldsymbol{\gamma}\).
We take \(\sigma_n>0\) whenever \(\alpha_n>0\); when \(\alpha_n=0\), the corresponding kernel term is defined to be zero.
The damped tangent Jacobian satisfies
\begin{equation}
\label{e:pk_tangent}
\begin{split}
\delta x_{n+1}^{\boldsymbol{\gamma}}
=(1-\alpha_n\Dt)\delta x_n^{\boldsymbol{\gamma}}
+\nabla F_n\,\delta x_n^{\boldsymbol{\gamma}}\Dt
+\delta F_n^{\boldsymbol{\gamma}}\Dt
+\DB_n\left((\delta x_n^{\boldsymbol{\gamma}})^T\nabla\sigma_n
+\delta\sigma_n^{\boldsymbol{\gamma}}\right)^T.
\end{split}
\end{equation}
Each column of \(\delta x_n^{\boldsymbol{\gamma}}\) is the path-kernel tangent for one parameter direction.

For a smooth terminal observable \(\Psi\), let
\(\widetilde\Psi_N:=\Psi(x_N)-\E{\Psi(x_N)}\).
The exact finite-time path-kernel identity is
\begin{equation}
\label{e:pk_identity}
\delta \E{\Psi(x_N^{\boldsymbol{\gamma}})}
=\E{(\delta x_N^{\boldsymbol{\gamma}})^T\nabla\Psi(x_N)
+\widetilde\Psi_N\sum_{n=0}^{N-1}
\frac{\alpha_n}{\sigma_n}(\delta x_n^{\boldsymbol{\gamma}})^T\DB_n}.
\end{equation}
The first term is the pathwise derivative.
The second transfers the part \(\alpha_n\delta x_n^{\boldsymbol{\gamma}}\) removed from the tangent dynamics to the probability kernel.
Thus \(\alpha_n>0\) damps unstable tangent directions without changing the exact parameter gradient.

\subsubsection{Gaussian initial-law extension}
We replace a deterministic initial condition by
\begin{equation*}
x_{-1}^{\boldsymbol{\gamma}}
:=\mu_{\mathrm{ic}}^{\boldsymbol{\gamma}},
\qquad
x_0^{\boldsymbol{\gamma}}
=x_{-1}^{\boldsymbol{\gamma}}
+\sigma_{\mathrm{ic}}^{\boldsymbol{\gamma}}B_{\mathrm{ic}},
\qquad
B_{\mathrm{ic}}\sim\cN(0,I),
\end{equation*}
where \(B_{\mathrm{ic}}\) is independent of the subsequent model-noise increments.
At the baseline parameter value,
\(x_0\sim\cN(\mu_{\mathrm{ic}},\sigma_{\mathrm{ic}}^2I)\).
We regard this draw as one auxiliary transition from time \(-1\) to time \(0\) and assign it a dimensionless schedule \(\alpha_{\mathrm{ic}}\in[0,1]\).
We take \(\sigma_{\mathrm{ic}}>0\) whenever \(\alpha_{\mathrm{ic}}>0\), and define the initial kernel term to be zero when \(\alpha_{\mathrm{ic}}=0\).
The initial tangent is
\[
\delta x_0^{\boldsymbol{\gamma}}
=(1-\alpha_{\mathrm{ic}})
\delta\mu_{\mathrm{ic}}^{\boldsymbol{\gamma}}
+B_{\mathrm{ic}}
(\delta\sigma_{\mathrm{ic}}^{\boldsymbol{\gamma}})^T.
\]
There is no factor \(\Dt\) because this is one auxiliary transition rather than one physical Euler step.

Applying the preceding tangent recursion for \(n=0,\ldots,N-1\) gives
\begin{equation*}
\begin{split}
\delta\E{\Psi(x_N^{\boldsymbol{\gamma}})}
=\mathbb{E}\bigg[(\delta x_N^{\boldsymbol{\gamma}})^T\nabla\Psi(x_N)
+\widetilde\Psi_N\left(
\frac{\alpha_{\mathrm{ic}}}{\sigma_{\mathrm{ic}}}
(\delta\mu_{\mathrm{ic}}^{\boldsymbol{\gamma}})^TB_{\mathrm{ic}}
+\sum_{n=0}^{N-1}\frac{\alpha_n}{\sigma_n}
(\delta x_n^{\boldsymbol{\gamma}})^T\DB_n\right)\bigg].
\end{split}
\end{equation*}
The choice \(\alpha_{\mathrm{ic}}=0\) gives the usual reparameterization derivative, whereas \(\alpha_{\mathrm{ic}}=1\) transfers the mean derivative completely to the Gaussian kernel.

\subsubsection{Full-window APK adjoint}
The branch-discovery objective depends on the whole sampled path rather than only on \(x_N\).
On the assimilation window, each APK sample path evolves according to
\begin{equation*}
x_{n+1}^{\boldsymbol{\gamma}}
=x_n^{\boldsymbol{\gamma}}+F^{\boldsymbol{\gamma}}(x_n^{\boldsymbol{\gamma}})\Dt
+\sigma^{\boldsymbol{\gamma}}\DB_n
+\xi^{\boldsymbol{\gamma}}(x_n^{\boldsymbol{\gamma}},x_n^{*,\boldsymbol{\gamma}})\Dt.
\end{equation*}
Here \(\xi^{\boldsymbol{\gamma}}\) is the control applied to the sampled dynamics; this paper uses the linear feedback form introduced in the next subsection.
The parameter vector may contain parameters in \(F^{\boldsymbol{\gamma}}\), \(\sigma^{\boldsymbol{\gamma}}\), and \(\xi^{\boldsymbol{\gamma}}\), together with \(\mu_{\mathrm{ic}}\), \(\sigma_{\mathrm{ic}}\), and the auxiliary path \(x_{0:N-1}^{*,\boldsymbol{\gamma}}\).
The model starts from the Gaussian state above, so \(\sigma_{\mathrm{ic}}\) is distinct from the model diffusion scale \(\sigma\).

The full-window loss is
\begin{equation}
\label{e:apk_da_loss}
\begin{split}
\Phi^{\boldsymbol{\gamma}}(x_{0:N}^{\boldsymbol{\gamma}})
:=\frac{1}{2T}\sum_{n=0}^{N-1}\bigg(|\phi(x_n^{\boldsymbol{\gamma}})-y_n|^2
+C|\xi^{\boldsymbol{\gamma}}(x_n^{\boldsymbol{\gamma}},x_n^{*,\boldsymbol{\gamma}})|^2\bigg)\Dt,
\qquad T=N\Dt.
\end{split}
\end{equation}
The constant \(C>0\) balances observation mismatch and dynamical correction.
An initial value of \(C\) may be chosen by balancing these two terms using the scaling formula in \cite{apk}.
The required averages can be estimated from an independent simulation of the base primal system.

At \(\boldsymbol{\gamma}=0\), set
\(\xi_n:=\xi(x_n,x_n^*)\),
\(\widetilde\Phi:=\Phi-\E{\Phi}\), and
\(\delta\xi_n^{\boldsymbol{\gamma}}:=\delta[\xi^{\boldsymbol{\gamma}}(x_n,x_n^{*,\boldsymbol{\gamma}})]\), where \(x_n\) is held fixed in the last derivative.
Applying \Cref{e:pk_tangent,e:pk_identity} to the full-window loss gives its tangent path-kernel formula.
The adjoint formula below then follows by the same tangent-adjoint duality argument as in \cite{apk}.

Define
\[
M_n:=(1-\alpha_n\Dt)I+\nabla F_n\Dt+\nabla_x\xi_n\Dt
\]
and
\[
\omega_n:=\nabla\phi_n^T(\phi(x_n)-y_n)\Dt
+\frac{T\alpha_n}{\sigma}\widetilde\Phi\,\DB_n
+C\nabla_x\xi_n^T\xi_n\Dt.
\]
The APK adjoint is the backward recursion
\[
\nu_N=0,
\qquad \nu_n=M_n^T\nu_{n+1}+\omega_n,
\qquad n=N-1,\ldots,0.
\]
The auxiliary initial step is
\[
\nu_{-1}=(1-\alpha_{\mathrm{ic}})\nu_0
+\frac{T\alpha_{\mathrm{ic}}}{\sigma_{\mathrm{ic}}}
\widetilde\Phi B_{\mathrm{ic}}.
\]

\begin{formula}[APK gradients for full-window data assimilation]
\label{t:apk_da}
For the objective in \Cref{e:apk_da_loss}, tangent-adjoint duality gives
\begin{equation*}
\delta\E{\Phi^{\boldsymbol{\gamma}}}
=\frac1T\mathbb{E}\left[
\begin{aligned}
&(\delta\mu_{\mathrm{ic}}^{\boldsymbol{\gamma}})^T\nu_{-1}
+\delta\sigma_{\mathrm{ic}}^{\boldsymbol{\gamma}}(B_{\mathrm{ic}}\cdot\nu_0)
+\sum_{n=0}^{N-1}(\delta F_n^{\boldsymbol{\gamma}})^T\nu_{n+1}\Dt
\\
&+\delta\sigma^{\boldsymbol{\gamma}}\sum_{n=0}^{N-1}\DB_n\cdot\nu_{n+1}
+\sum_{n=0}^{N-1}(\delta\xi_n^{\boldsymbol{\gamma}})^T
(\nu_{n+1}+C\xi_n)\Dt
\end{aligned}
\right].
\end{equation*}
\end{formula}

We can read off the gradient directions from the preceding expression.
For example, the gradients with respect to \(\mu_{\mathrm{ic}}\) and \(\sigma_{\mathrm{ic}}\) are
\[
\left[\delta\E{\Phi^{\boldsymbol{\gamma}}}\right]_{\mu_{\mathrm{ic}}}
=\frac1T\E{\nu_{-1}},
\qquad
\left[\delta\E{\Phi^{\boldsymbol{\gamma}}}\right]_{\sigma_{\mathrm{ic}}}
=\frac1T\E{B_{\mathrm{ic}}\cdot\nu_0}.
\]
All expectations include the initial draw and all subsequent model-noise increments.
In computation, these expectations are estimated using \(L\) stochastic sample paths for each APK run at each optimization update.
All \(L\) paths use the same current parameters.
In particular, let \(\Phi_l\) be the loss of path \(l\), and define
\begin{equation}
\label{e:apk_batch_centering}
\bar\Phi:=\frac1L\sum_{s=1}^L\Phi_s,
\qquad
\widehat{\widetilde\Phi}_l:=\Phi_l-\bar\Phi.
\end{equation}
Because \(\bar\Phi\) includes path \(l\), \(\widehat{\widetilde\Phi}_l\) is not an unbiased estimator of \(\widetilde\Phi_l\).
However, it has the lowest variance among the standard same-batch centering choices and is more stable for small \(L\).
We substitute \(\widehat{\widetilde\Phi}_l\) for \(\widetilde\Phi\) in the kernel terms and average the resulting gradients over the \(L\) paths.
When the kernel terms vanish, one sample path is sufficient.

\subsection{Exploration and exploitation under prescribed schedules}
\label{s:apk_exploration}
\label{s:apk_linear_correction}

The correction used in this paper is
\begin{equation}
\label{e:linear_correction}
\xi_n
=
\gamma'(x_n^*-x_n).
\end{equation}
The auxiliary path \(x^*\) is the center of this linear attraction.
The path \(x_{0:N-1}^*\) and initial mean \(\mu_{\mathrm{ic}}\) are optimized with the APK gradient, while \(\gamma'\), \(\sigma\), and \(\sigma_{\mathrm{ic}}\) follow prescribed values.
The correction is penalized by the fixed coefficient \(C\) in \Cref{e:apk_da_loss}.

At optimization update \(k\), the reported method uses
\begin{equation}
\label{e:prescribed_parameter_schedule}
\gamma'^{(k)}=\gamma',
\qquad
\alpha^{(k)}=\sigma^{(k)}=\sigma_{\mathrm{ic}}^{(k)}
=\alpha_0\max\left\{1-\frac{k}{K},0\right\},
\qquad
\alpha_{\mathrm{ic}}^{(k)}
=\max\left\{1-\frac{k}{K},0\right\}.
\end{equation}
Here \(\alpha_0>0\) sets the initial path-kernel damping and stochastic scales, and \(K\) is the number of exploration updates.
Thus the path-kernel damping and both stochastic scales decay together, while the linear correction remains fixed.

For the linear correction, \(\nabla_x\xi_n=-\gamma'I\), so the state-propagation matrix in the adjoint at update \(k\) is
\begin{equation*}
M_n^{(k)}
=I+\Dt\left(
\nabla F_n-\left[\alpha^{(k)}+\gamma'\right]I
\right).
\end{equation*}
The pathwise damping therefore decreases from \(\alpha_0+\gamma'\) to the fixed correction scale \(\gamma'\), rather than remaining constant.

At initialization,
\[
\alpha=\sigma=\sigma_{\mathrm{ic}}=\alpha_0,
\qquad
\alpha_{\mathrm{ic}}=1.
\]
The large \(\alpha\) stabilizes the long-window tangent and adjoint by removing the unstable part of the pathwise sensitivity and transferring it, through the path-kernel identity, to the centered probability-kernel terms
\[
\frac{T\alpha}{\sigma}\widetilde\Phi\,\DB_n,
\qquad
\frac{T\alpha_{\mathrm{ic}}}{\sigma_{\mathrm{ic}}}
\widetilde\Phi B_{\mathrm{ic}}.
\]
At the same time, the large values of \(\sigma\) and \(\sigma_{\mathrm{ic}}\) generate broad sample paths, which may travel far from the center \(x^*\).
The centered kernel contribution makes this exploration directional \cite{apk}: if a stochastic path has below-average full-window loss, its realized noise increments update the center and the other optimization variables so that similar excursions become more likely.
Thus a useful random path can move the optimization toward a distant low-loss branch before the center path has reached that branch through local updates.

As the prescribed \(\alpha\), \(\sigma\), and \(\sigma_{\mathrm{ic}}\) decrease, stability is increasingly supplied by the fixed deterministic attraction toward \(x^*\), rather than by transferring sensitivity to the probability kernel.
In the late regime, \(\gamma'\) remains fixed while \(\alpha\), \(\sigma\), and \(\sigma_{\mathrm{ic}}\) vanish.
The sampled paths concentrate around the center and their centered loss differences become small, so random excursions receive little additional reward from the kernel terms.
The optimization is then dominated by pathwise adjoint updates of \(x^*\) and \(\mu_{\mathrm{ic}}\): it mainly rewards a lower-loss change in the center of the path.
This is the exploitation regime.
Although the ratio \(\alpha/\sigma\) need not itself vanish, the kernel contribution vanishes in the deterministic limit because the sampled paths and their centered loss fluctuations collapse.

The fixed correction strength can be related to local instability.
For a small error \(e_n\) from the auxiliary path, the corrected error dynamics is locally
\[
e_{n+1}
\approx
\left[I+\Dt\bigl(\nabla F(x_n)-\gamma'I\bigr)\right]e_n.
\]
Its leading transverse growth rate is approximately \(\lambda_{\max}-\gamma'\), where \(\lambda_{\max}\) is the largest Lyapunov exponent of the uncorrected dynamics.
Thus choosing \(\gamma'\) comfortably above \(\lambda_{\max}\) makes the corrected error dynamics contractive, although the underlying system remains chaotic.
The initial scale \(\alpha_0\) should be large enough to provide early gradient stabilization and stochastic exploration, while \(K\) determines how long this regime lasts.
The particular values used in the experiments are given in \Cref{s:numeric_setting}.
The early kernel-differentiation regime helps move \(x^*\) and \(\mu_{\mathrm{ic}}\) toward a useful branch from a rough initial guess, and the fixed correction-stabilized path-perturbation regime helps concentrate broad stochastic exploration into an accurate path.

\subsection{Population mixing: construct mixed path proposals}
\label{s:theory_population}
\label{s:apk_population_mix}

Population mixing acts across independently optimized full-window paths.
Population methods have previously been used to address local minima in variational data assimilation.
Ahrens \cite{Ahrens1999} represents a weak-constraint trajectory by a chromosome, evaluates each candidate using one full-window objective, and randomly exchanges control variables through uniform crossover.
Although global fitness determines which candidates are more likely to survive, the crossover itself does not identify which portions of each trajectory are successful.
Our method instead assigns credit at each physical time and state coordinate.
Temporally averaged coordinatewise losses determine how the locally successful portions of independently optimized paths are combined.
Subsequent APK updates restore the dynamical coherence of the mixed proposal.
Thus, the essential mechanism is informed spatiotemporal recombination rather than random crossover between globally ranked trajectories.
Note that the population runs used here are distinct from the stochastic sample paths within one APK update.
The stochastic paths are averaged at every optimization update to estimate the expected gradient for one run, whereas the population runs explore different basins and are mixed only at prescribed intervals, every \(1000\) updates in our experiments.

The population is initialized from \(R_{\rm apk}\) independent model rollouts.
For run \(r\), each coordinate of the initial mean is drawn from the following law, where \(F_{\rm forc}\) denotes the scalar forcing:
\[
\mu_{{\rm ic},j}^{(r)}\sim\mathcal N(F_{\rm forc}/2,2^2),
\qquad j=1,\ldots,M,
\]
and the auxiliary path starts from this mean and follows the deterministic model:
\[
x_0^{*,(r)}=\mu_{\mathrm{ic}}^{(r)},
\qquad
x_{n+1}^{*,(r)}=x_n^{*,(r)}+F(x_n^{*,(r)})\Dt.
\]
After the initial auxiliary path \(x^{*,(r)}\) is constructed for each run \(r\), the runs are optimized independently between prescribed population-mixing updates.

Local losses assign credit to every run at each physical time and state coordinate.
At each prescribed mixing update, we first evaluate the deterministic representative path of every run.
This path is obtained by setting \(B_{\mathrm{ic}}=0\) and \(\DB_n=0\) in the controlled dynamics, so that \(x_0^{(r)}=\mu_{\mathrm{ic}}^{(r)}\), while retaining the correction.
For \(r=1,\ldots,R_{\rm apk}\), the deterministic representative path and correction are written as
\[
x_{0:N}^{(r)},
\qquad
\xi_n^{(r)}:=\gamma'^{(r)}(x_n^{*,(r)}-x_n^{(r)}).
\]
For time \(n\) and state coordinate \(j\), define the local spatiotemporal loss as
\begin{equation*}
\ell_{n,j}^{(r)}
:=
\frac12\left(
\left|\phi(x_n^{(r)})-y_n\right|^2
+CM\left|\xi_{n,j}^{(r)}\right|^2
\right).
\end{equation*}
The observation discrepancy is shared by all state coordinates at the same time, whereas the correction discrepancy is coordinatewise.
Averaging over the \(M\) state coordinates gives
\[
\frac1M\sum_{j=1}^M\ell_{n,j}^{(r)}
=
\frac12\left(
\left|\phi(x_n^{(r)})-y_n\right|^2
+C\left|\xi_n^{(r)}\right|^2
\right),
\]
the time-\(n\) integrand of the APK objective.

The mixing criterion first averages each coordinate loss over a Gaussian temporal neighborhood:
\begin{equation*}
\overline\ell_{n_0,j}^{(r)}
:=
\frac{
\displaystyle\sum_{n=0}^{N-1}
\exp{-\frac{(t_n-t_{n_0})^2}{S^2}}
\ell_{n,j}^{(r)}
}{
\displaystyle\sum_{n=0}^{N-1}
\exp{-\frac{(t_n-t_{n_0})^2}{S^2}}
}.
\end{equation*}
The bell radius \(S\) is shared with restart selection below and is set by the local dynamical time scale, estimated from an independent simulation of the base primal system.
Let \(\tau_{\mathrm{dec}}\) be the first lag at which the mean coordinate autocorrelation along a representative model trajectory reaches \(e^{-1}\):
\begin{equation*}
\frac1M\sum_{j=1}^M
\operatorname{Corr}\left(x_t^j,x_{t+\tau_{\mathrm{dec}}}^j\right)
=e^{-1}.
\end{equation*}
We set \(S=\tau_{\mathrm{dec}}/2\), which smooths stepwise loss fluctuations while keeping the average local relative to the dynamics.

The mixing criterion then converts these temporally averaged losses into soft-min weights across runs:
\begin{equation*}
w_{n,j}^{(r)}
:=
\frac{
\exp{-\beta\left(
\overline\ell_{n,j}^{(r)}
-\min_q\overline\ell_{n,j}^{(q)}
\right)}
}{
\displaystyle\sum_{s=1}^{R_{\rm apk}}
\exp{-\beta\left(
\overline\ell_{n,j}^{(s)}
-\min_q\overline\ell_{n,j}^{(q)}
\right)}
}.
\end{equation*}
Here \(\beta>0\) is the inverse temperature and controls how strongly the proposal favors the lowest local losses.
Subtracting the pointwise minimum prevents numerical underflow without changing the normalized weights.

The population-wise mixed proposal for the auxiliary path and initial mean is determined by
\begin{equation*}
\widetilde x_{n,j}^*
=
\sum_{r=1}^{R_{\rm apk}}w_{n,j}^{(r)}x_{n,j}^{*,(r)},
\qquad
\widetilde\mu_{{\rm ic},j}
=
\sum_{r=1}^{R_{\rm apk}}w_{0,j}^{(r)}
\mu_{{\rm ic},j}^{(r)}.
\end{equation*}
The weights vary in both time and space, so one run need not supply the entire proposal.
The reported method does not mix scalar parameters: \(\gamma'\) remains fixed, and \(\alpha\), \(\sigma\), and \(\sigma_{\mathrm{ic}}\) retain their prescribed values at the current optimization update.

Each mixing update replaces only the current worst population run.
Specifically, the run with the largest current deterministic full-window objective receives the mixed proposal; all other runs remain unchanged.
Therefore an existing best run is retained, and the smallest objective already present in the population cannot increase merely because mixing is performed.
The mixed proposal itself need not immediately reduce the full objective, because recombination may introduce mismatch between neighboring times and coordinates; subsequent APK updates repair this mismatch.
After replacement, the replaced run keeps its own future randomness.
Only runs with finite states, parameters, and local losses contribute to the mixture.
This operation is inherently noncausal.
It evaluates complete paths using the entire assimilation window and revises different times and coordinates retrospectively.

\subsection{Loss-selected EnKF restart and continuation}
\label{s:local_filtering_after_apk}
\label{s:enkf_after_init}

The optional EnKF continuation separates global branch discovery from local tracking.
Before enough observations have accumulated, the causal filtering marginal can contain several plausible branches, and an EnKF initialized on one of them may remain on the wrong branch or diverge.
APK instead uses the offline assimilation window retrospectively to identify a low-loss branch and restart state.
If the selected branch is locally dominant and concentrated, approximating its uncertainty by a single-mode model such as a Gaussian is typically reasonable, and the EnKF only needs to preserve a branch that has already been identified.

Final run selection is based on the deterministic full-window objective.
After the final optimization update, we select the run with the smallest best attained value of this objective, and denote its selected path by \(x_{0:N}^{\mathrm{apk}}\).

The restart-time criterion begins with the instantaneous loss along the selected APK path:
\[
\ell_n^{\mathrm{apk}}
:=
\frac12\left(
|\phi(x_n^{\mathrm{apk}})-y_n|^2
+C|\xi_n^{\mathrm{apk}}|^2
\right),
\qquad 0\leq n<N,
\]
where \(\xi_n^{\mathrm{apk}}\) is the optimized correction.
The full-window objective is the temporal average of this instantaneous loss, so a small objective guarantees a small instantaneous loss somewhere in the window.

Temporal averaging prevents an isolated loss fluctuation from determining the restart time.
Using the same bell-shaped average and radius \(S\) as in population mixing, for \(0\leq n_0\leq N\), set
\begin{equation*}
\overline\ell_{n_0}^{\mathrm{apk}}
:=
\frac{
\displaystyle\sum_{n=0}^{N-1}
\exp{-\frac{(t_n-t_{n_0})^2}{S^2}}
\ell_n^{\mathrm{apk}}
}{
\displaystyle\sum_{n=0}^{N-1}
\exp{-\frac{(t_n-t_{n_0})^2}{S^2}}
}.
\end{equation*}

Restart selection is restricted to the interior because the loss is less constrained near the two endpoints of the reconstruction window.
Rather than add an endpoint weight to the bell-averaged loss, we restrict the candidate restart times to be at least one model-time unit from both endpoints:
\begin{equation}
\label{e:apk_restart_loss}
n_*
:=
\mathop{\mathrm{arg\,min}}_{\substack{0\leq n_0\leq N\\
1\leq t_{n_0}\leq t_N-1}}
\overline\ell_{n_0}^{\mathrm{apk}}.
\end{equation}
We write \(t_*:=t_{n_*}\) for the corresponding restart time.

The resulting restart score is computable without the true path, but it does not automatically determine the full-state error under partial observations.
Its interpretation requires the selected branch to be identifiable from the observations and dynamical consistency over this temporal neighborhood.
When this condition holds, a small selected score \(\overline\ell_{n_*}^{\mathrm{apk}}\) indicates a small state error at the restart time.

The reconstructed full-window path is the primary APK output.
If the restart state selected by \Cref{e:apk_restart_loss} lies on a locally viable branch, an optional EnKF can track the selected branch forward from that state.
Combining variational and ensemble assimilation is already standard in hybrid 4D-Var/EnKF methods \cite{ZhangZhang2012E4DVar}.
Our contribution here is a truth-free local loss criterion for selecting the EnKF restart time within the offline assimilation window.

The EnKF ensemble is initialized with a small Gaussian spread around the selected APK state, replayed from \(n_*\) through the remaining observations in the offline assimilation window, and then continued causally beyond \(T\).
The spread is an online modeling choice; we use a small isotropic spread in the numerical experiments.

The EnKF then performs inexpensive local tracking using ensemble cross-covariances.
These covariances can transfer observed innovations to hidden coordinates, but they do not select among globally observation-compatible branches.
Continuation is therefore expected to remain reliable only while the ensemble stays concentrated within the selected branch and the dynamics and observation map remain locally well approximated.
In the numerical comparison, we use the same locally stable inflation for the APK-initialized and rough-prior EnKFs; its calibration is reported in Section~\ref{s:numeric_setting}.
A large innovation relative to the predicted observation spread can indicate that the EnKF has left this local regime, in which case APK can be rerun on a recent observation window.
If online continuation is not viable, the reconstructed APK path remains available.

\subsection{Complete procedure}
\label{s:algorithm_details}

The complete procedure has three stages: population APK optimization with occasional mixing, loss-based selection of a path and restart time, and optional EnKF continuation.
\Cref{a:apk_enkf} summarizes these stages; the numerical settings are given in Section~\ref{s:numeric_setting}.

\begin{algorithm}[!htbp]
\caption{Full-window branch discovery with optional EnKF continuation}
\label{a:apk_enkf}
\begin{algorithmic}[1]
\Require Observations \(y_{0:N}\) on the offline assimilation window, model \(F\), observation map \(\phi\), APK and mixing settings, and optionally online observations and EnKF settings.
\Ensure Selected full-window path \(x_{0:N}^{\mathrm{apk}}\), restart index \(n_*\), and optionally an online EnKF ensemble.
\State Independently initialize each APK run by drawing its initial mean and generating its auxiliary path by a deterministic model rollout.
\For{each APK optimization update}
  \For{each population run}
    \State Draw the stochastic sample paths and apply one APK mini-batch update using \Cref{t:apk_da,e:apk_batch_centering,e:linear_correction,e:prescribed_parameter_schedule}.
  \EndFor
  \If{this is a prescribed mixing update}
    \State Form the spatiotemporal mixed proposal and use it to replace the current worst run, as in \Cref{s:apk_population_mix}.
  \EndIf
\EndFor
\State Select the run with the smallest best attained deterministic objective, evaluated with the common penalty \(C\), and denote its optimized path by \(x_{0:N}^{\mathrm{apk}}\).
\State Select \(n_*\) by minimizing the local loss in \Cref{e:apk_restart_loss} over times at least one model-time unit from both endpoints.
\If{an online continuation is used}
  \State Initialize an ensemble with a small Gaussian spread around the APK restart state.
  \State Replay the EnKF through the remaining observations in the offline assimilation window.
  \State Continue the EnKF causally beyond \(T\).
\EndIf
\end{algorithmic}
\end{algorithm}

\section{Numerical Results}
\label{s:numeric}

We use two \(40\)-dimensional Lorenz--96 settings.
The main linear-observation experiment compares offline APK reconstruction with population weak-\(\mathrm{4D\text{-}Var}_x\) and compares an EnKF initialized from the APK reconstruction with the same EnKF started from a rough prior.
A separate squared-observation experiment tests whether the same conclusions persist for a nonlinear observation map.
All statistics are computed from twenty independently seeded experiments in each setting.
For the main figures, we use the experiment whose APK path RMSE \(0.779\) is closest to the twenty-experiment median \(0.800\).

\subsection{Shared Lorenz--96 setting and evaluation}
\label{s:numeric_setting}

We use the Euler-discretized Lorenz--96 dynamics
\begin{equation}
\label{e:l96_numeric}
x_{j,n+1}
=
x_{j,n}
+\Dt\left[
(x_{j+1,n}-x_{j-2,n})x_{j-1,n}
-x_{j,n}+F_{\rm forc}
\right],
\end{equation}
with cyclic indices \(j=0,\ldots,39\), forcing \(F_{\rm forc}=8\), and \(\Dt=0.005\).
The true state is integrated for \(4000\) steps before the offline assimilation window.
Both the data-generating and inference models are deterministic, use the same known forcing, and cover \(0\leq t\leq 2T=10\).

We observe one coordinate and then miss four, repeated around the cyclic state:
\begin{equation}
\label{e:l96_obs_indices}
\mathcal I_{\rm obs}
=
\{0,5,10,\ldots,35\}.
\end{equation}
Thus \(8\) coordinates are observed and \(32\) are hidden.
Observations are available at every Euler step and contain independent Gaussian noise,
\[
y_n=P_{\rm obs}x_n+0.3\eta_n,
\qquad
\eta_n\sim\mathcal N(0,I_8).
\]
Both EnKF updates use the correct observation covariance \(0.3^2I_8\).

APK and weak-\(\mathrm{4D\text{-}Var}_x\) use the offline assimilation window \(0\leq t\leq T=5\).
Each method has \(16\) independently initialized population members.
APK performs \(5000\) stochastic mini-batch updates, while weak-\(\mathrm{4D\text{-}Var}_x\) receives \(20000\) deterministic optimization updates.
For every member, the initial state is drawn coordinatewise from
\[
\mathcal N(F_{\rm forc}/2,2^2)
\]
and propagated with the deterministic model to form the initial auxiliary path.
Both methods update each state coordinate using an objective-scaled clipping rule that limits its predicted first-order objective change to \(1/40\) of the current objective; see \cite{apk} for details.
The uncapped base step size is \(0.5\).
APK applies the resulting step directly, while weak-\(\mathrm{4D\text{-}Var}_x\) uses Armijo backtracking and leaves the path unchanged when no finite acceptable trial is found.
For APK, the initial-mean step uses three times this base value, and the auxiliary-path step includes the \(1/\Dt\) compensation for its \(O(\Dt)\) gradient before clipping.
The per-update change of each initial-mean coordinate is bounded by \(0.3\).
We fix \(\gamma'=4\) and prescribe
\[
\alpha^{(k)}=\sigma^{(k)}=\sigma_{\mathrm{ic}}^{(k)}
=4\max\left\{1-\frac{k}{3000},0\right\},
\qquad
\alpha_{\mathrm{ic}}^{(k)}
=\max\left\{1-\frac{k}{3000},0\right\}.
\]
Thus \(\gamma'=4\) remains fixed, while \(\alpha\), \(\sigma\), and \(\sigma_{\mathrm{ic}}\) decrease together and reach zero at update \(3000\).
Mixing begins at update \(2010\), repeats every \(1000\) updates, uses inverse temperature \(40\), and uses a Gaussian temporal bell of radius
\[
S=0.135,
\]
one half of the measured \(0.270\) e-folding decorrelation time.
APK uses two stochastic sample paths while the kernel terms are active and one path after update \(3000\).
The fixed correction penalty is
\[
C\approx0.00716.
\]

Writing \(F_{\rm L96}(x)\) for the bracketed Lorenz--96 drift in \Cref{e:l96_numeric}, the weak-\(\mathrm{4D\text{-}Var}_x\) comparator directly optimizes every state in the path using the objective
\[
\frac1{2N}\sum_{n=0}^{N-1}
\left[
\frac{|P_{\rm obs}x_n-y_n|^2}{\sigma_o^2}
+
\frac{|x_{n+1}-x_n-\Dt F_{\rm L96}(x_n)|^2}{\sigma_q^2}
\right],
\qquad
\sigma_o=0.3,
\qquad
\sigma_q=\frac{\sigma_o\Dt}{\sqrt C}.
\]
This choice gives the observation discrepancy and model residual the same relative weighting as the APK observation and correction terms.
The weak-\(\mathrm{4D\text{-}Var}_x\) comparator also receives the same population-mixing schedule.
Thus, it is not disadvantaged in either the number of optimization updates or the opportunities for population mixing.

For either reconstructed path, the restart time minimizes the bell-averaged local loss subject to
\[
1\leq t_*\leq T-1.
\]
The truth is not used to choose the population member or \(t_*\).
For the online comparison, both EnKFs have \(120\) members, inflation
\[
\rho=1.01,
\]
no process noise, and no forcing uncertainty.
We screened \(\rho\) using accurately initialized EnKFs to identify a value that permits nonexplosive local tracking, and then fixed the same value for both online methods.
We also tested nine values \(0.98\leq\rho\leq1.04\) on ten independently seeded validation experiments.
For \(\rho-1=0,0.005,0.01,0.015\), every accurately initialized EnKF remains finite, with mean online RMSE between \(0.054\) and \(0.150\).
The corresponding rough-prior EnKFs are finite in \(5,1,8,0\) experiments, respectively; where defined, their mean RMSE over complete experiments is at least \(6.50\).
The accurately initialized EnKFs with the smaller inflation values \(0.98\) and \(0.99\) remain finite but have mean RMSE above \(5.2\).
Thus, several tested inflation values preserve accurate local tracking, but none makes the rough-prior EnKF accurate.
Among these locally stable values, \(\rho=1.01\) gives the largest number of finite rough-prior EnKF trajectories.
The APK-initialized ensemble has spread \(0.05\) around the APK state at \(t_*\), is replayed through the remainder of the offline assimilation window, and then continues online.
The ordinary EnKF starts at \(t=0\) from an independently generated climatological ensemble.

The true state \(x_n^{\rm true}\) is unavailable to every assimilation method and is used only to generate the synthetic observations and evaluate the reported errors.
Separate truth-initialized validation experiments are used only to select \(\rho\), not as inputs to the reported assimilation experiments.
At each physical time, the all-coordinate RMSE is
\[
\operatorname{RMSE}(t_n)
=
\left[
\frac1{40}
\left|\widehat x_n-x_n^{\rm true}\right|^2
\right]^{1/2}.
\]
The reported offline path RMSE is the root mean square of this quantity over \(1\leq t\leq T-1\), excluding the less-constrained endpoints.
The online RMSE is its root mean square over \(T\leq t\leq2T\).
The restart RMSE is evaluated at the selected \(t_*\).
We report the mean and sample standard deviation across experiments.
An experiment contributes to a method's aggregate whenever that method is finite throughout the stated interval; a nonfinite ordinary-EnKF trajectory is counted as a failure rather than included in its RMSE average.

We now report the computational cost.
The offline APK search carries most of the computational cost, while the subsequent EnKF continuation is relatively inexpensive.
The representative APK calculation ran its \(16\) population members concurrently on one UCI HPC3 node.
Ignoring the few population-mixing reevaluations, the APK stage used approximately \(128000\) loss-only stochastic forward paths, \(128000\) replayed stochastic forward paths with adjoint sweeps, and \(80000\) deterministic representative-path evaluations.
The APK-only wall time reported by the driver was approximately \(1.89\times10^3\) seconds, or \(31.6\) minutes; this excludes the subsequent weak-\(\mathrm{4D\text{-}Var}_x\) calculation.
After initialization, each EnKF observation step requires \(120\) model forecasts and one \(8\times8\) observation-space covariance solve.

\subsection{APK versus weak-4D-VarX}
\label{s:numeric_4dvarx}

The direct long-window implementations of strong-constraint \(\mathrm{4D\text{-}Var}\) and weak-\(\mathrm{4D\text{-}Var}_q\) do not produce finite comparators.
For both methods, the gradient becomes nonfinite at the first optimization step, before any reconstructed path is produced, so no trajectory or RMSE can be reported.
We therefore use weak-\(\mathrm{4D\text{-}Var}_x\) as the only variational comparator.

The finite offline comparison strongly favors APK over population weak-\(\mathrm{4D\text{-}Var}_x\).
Both methods remain finite in all twenty experiments.
APK is more accurate in every paired experiment, and its mean interior path RMSE is \(4.3\) times smaller, with a mean paired RMSE reduction of \(2.33\).
Nine of the twenty final APK paths have RMSE below \(0.5\), and thirteen are below \(1\).
Moreover, the computable local loss identifies substantially more accurate restart states along the APK paths than along the weak-\(\mathrm{4D\text{-}Var}_x\) paths, as summarized in \Cref{t:l96_offline_statistics}.

\begin{table}[!ht]
\centering
\begin{tabular}{lcccc}
\toprule
Quantity
& experiments
& mean \(\pm\) SD
& median
& range\\
\midrule
APK path
& \(20\)
& \(0.698\pm0.466\)
& \(0.800\)
& \([0.109,1.35]\)\\
Weak-\(\mathrm{4D\text{-}Var}_x\) path
& \(20\)
& \(3.02\pm0.419\)
& \(2.99\)
& \([2.37,4.06]\)\\
APK restart state
& \(20\)
& \(0.180\pm0.092\)
& \(0.153\)
& \([0.051,0.375]\)\\
Weak-\(\mathrm{4D\text{-}Var}_x\) restart state
& \(20\)
& \(2.56\pm0.844\)
& \(2.23\)
& \([1.29,4.57]\)\\
\bottomrule
\end{tabular}
\caption{
Offline all-coordinate RMSE over \(1\leq t\leq T-1\), and pointwise RMSE at the loss-selected restart time.
Every row uses all twenty independently seeded experiments.
}
\label{t:l96_offline_statistics}
\end{table}

To show a typical rather than the best case in the figures below, we select the experiment whose APK interior path RMSE, \(0.779\), is closest to the twenty-experiment median \(0.800\).

\begin{figure}[!ht]
\centering
\includegraphics[width=0.75\textwidth]{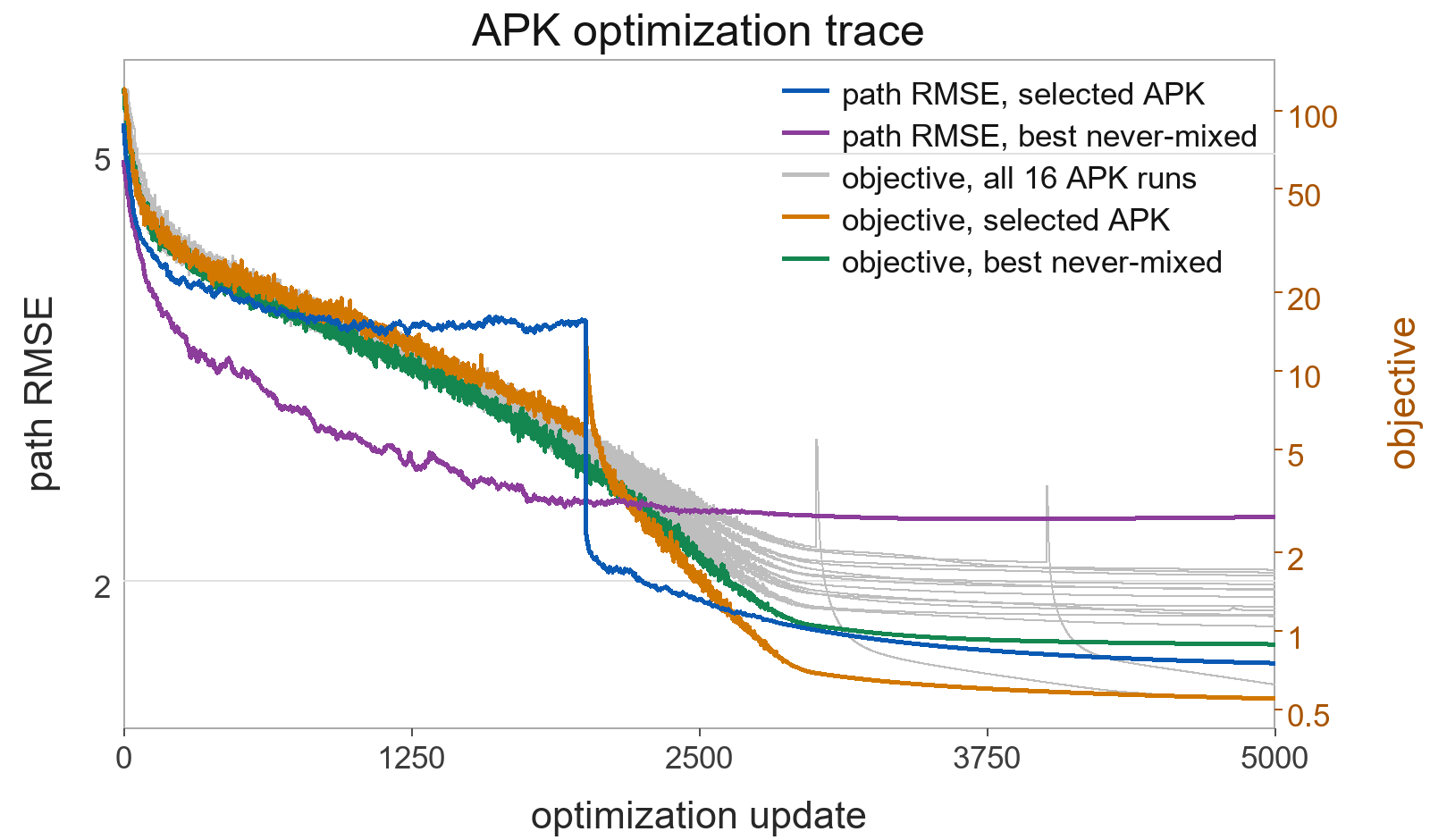}
\caption{
APK optimization in the representative experiment.
Gray curves show the objectives of all \(16\) runs, and the orange curve shows the objective-selected run.
The blue and purple curves show path RMSE for the selected and best never-mixed runs, respectively; the green curve shows the objective of the best never-mixed run.
The sharp decrease near update \(2010\) occurs when population mixing replaces the selected run.
}
\label{f:l96_apk_optimization_representative}
\end{figure}

The representative APK optimization illustrates how population mixing can improve a path.
The selected APK run is mixed once, at the first mixing event near update \(2010\), when its diagnostic path RMSE drops sharply in \Cref{f:l96_apk_optimization_representative}.
At the final update, this run has deterministic objective \(0.547\) and diagnostic whole-window path RMSE \(1.67\), compared with \(0.885\) and \(2.29\) for the best run that was never mixed.
Its evaluated interior path RMSE is \(0.779\), while the selected weak-\(\mathrm{4D\text{-}Var}_x\) path has interior RMSE \(2.75\).

\begin{figure}[!ht]
\centering
\includegraphics[width=0.75\textwidth]{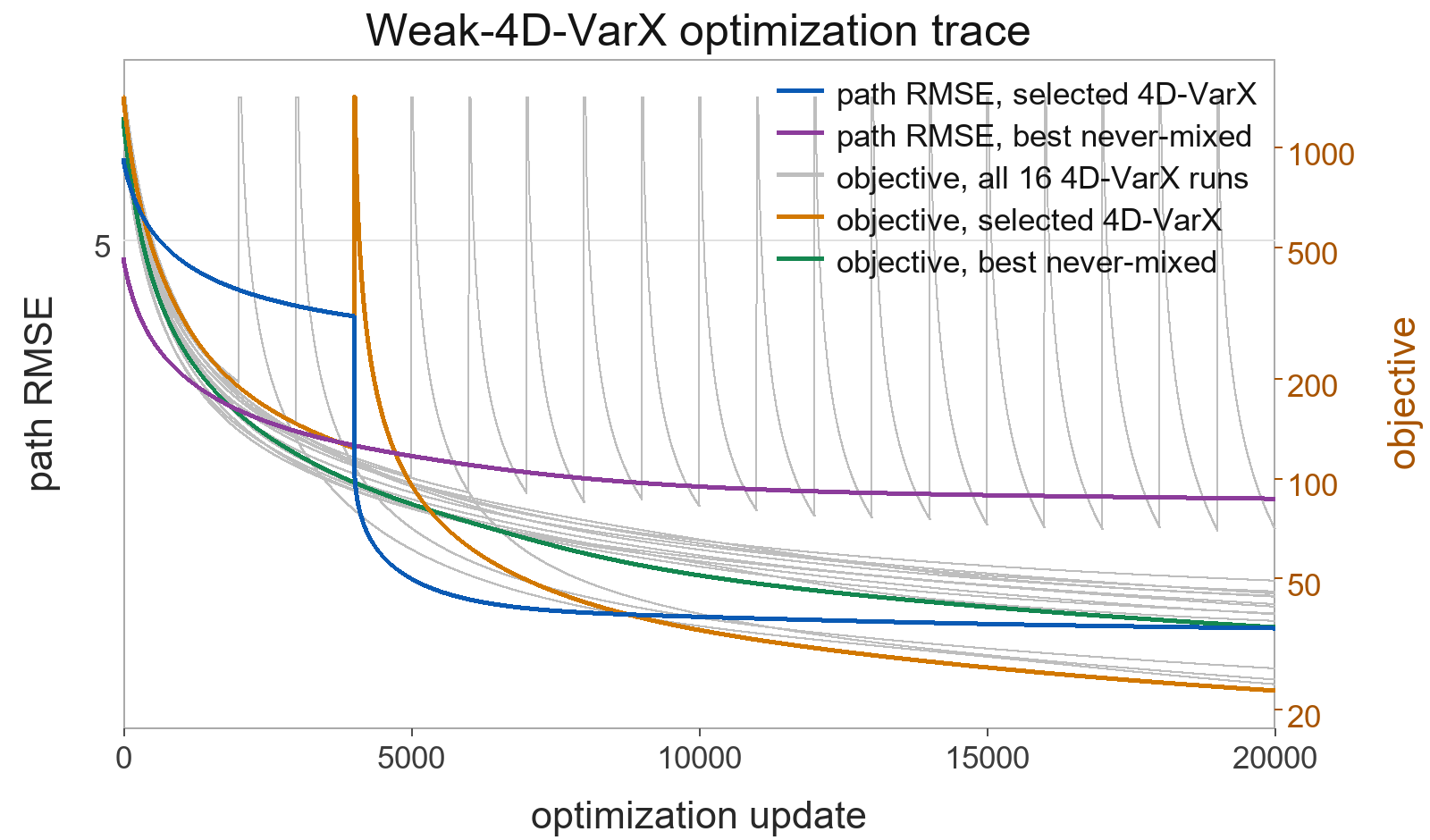}
\caption{
Weak-\(\mathrm{4D\text{-}Var}_x\) objective during population optimization in the representative experiment.
Line conventions follow \Cref{f:l96_apk_optimization_representative}.
The jumps occur when population mixing replaces a poor run.
}
\label{f:l96_weak_4dvarx_objective_trace}
\end{figure}

Weak-\(\mathrm{4D\text{-}Var}_x\) nevertheless achieves stable objective descent within the selected local basin.
In the representative experiment, its objective decreases from \(1.43\times10^3\) at the first reported update to \(22.8\) after \(20000\) updates, as shown in \Cref{f:l96_weak_4dvarx_objective_trace}.
Between population replacements, its Armijo updates produce stable descent.
The comparator therefore completes a long optimization with finite iterates, although this does not establish convergence to a global minimum.

\FloatBarrier
\subsection{APK--EnKF versus ordinary EnKF}
\label{s:numeric_enkf}

The online comparison strongly favors the APK-initialized EnKF over the ordinary EnKF.
The APK-initialized EnKF remains finite in all twenty experiments, whereas the ordinary EnKF becomes nonfinite in seven experiments.
Among the thirteen experiments in which both remain finite, the APK-initialized EnKF has a mean online RMSE \(63.8\) times smaller, with a mean absolute reduction of \(6.02\).
The complete statistics are given in \Cref{t:l96_filter_statistics}.

\begin{table}[!ht]
\centering
\begin{tabular}{lcccc}
\toprule
Method
& finite experiments
& mean \(\pm\) SD
& median
& range\\
\midrule
Ordinary EnKF
& \(13/20\)
& \(6.12\pm1.90\)
& \(6.08\)
& \([1.04,8.60]\)\\
APK--EnKF
& \(20/20\)
& \(0.0948\pm0.0050\)
& \(0.0947\)
& \([0.0844,0.104]\)\\
\bottomrule
\end{tabular}
\caption{
All-coordinate online RMSE over \(T\leq t\leq2T\).
The ordinary-EnKF statistics use its thirteen complete experiments; the other seven trajectories become nonfinite and are reported through the finite-experiment count.
}
\label{t:l96_filter_statistics}
\end{table}

\begin{figure}[!ht]
\centering
\includegraphics[width=0.8\textwidth]{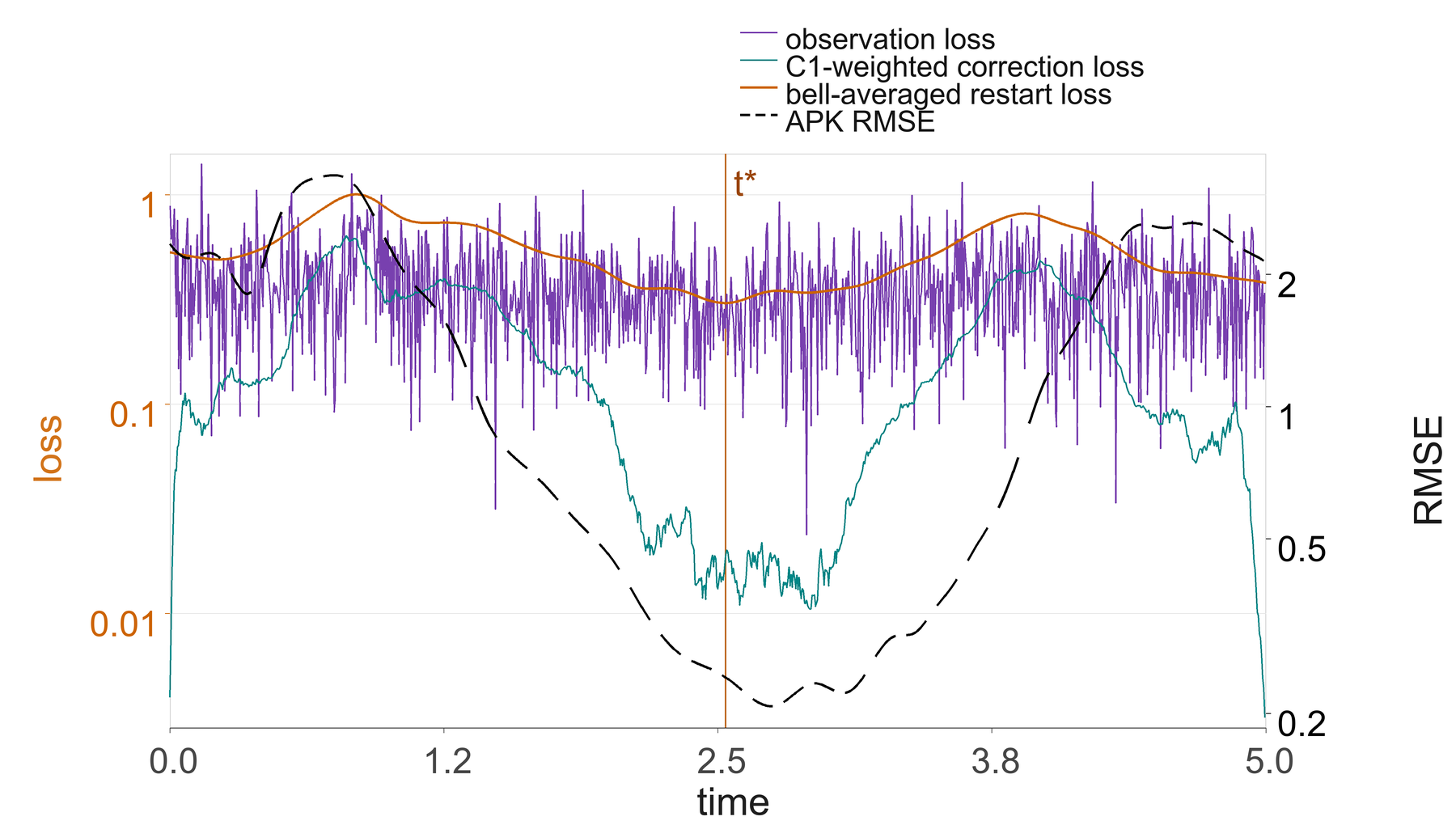}
\caption{
Restart-loss trajectory for the same representative APK reconstruction.
The purple and teal curves are the observation discrepancy and the \(C\)-weighted correction penalty, respectively.
The orange curve is their bell-averaged sum, which is minimized at the selected restart time \(t_*\approx2.54\).
The black dashed curve is the APK RMSE, shown on the right axis.
The three loss curves are computed from the observations and reconstructed APK path without using the true state; the RMSE is included only for evaluation.
}
\label{f:l96_restart_loss_representative}
\end{figure}

The representative experiment illustrates how an interior low-loss restart leads to accurate continuation.
Its bell-averaged loss selects the restart time \(t_*\approx2.54\), where the APK state RMSE is \(0.242\).
The APK--EnKF online RMSE is \(0.0945\).
The ordinary EnKF follows a high-error trajectory and eventually becomes nonfinite.
The black curve begins at \(t_*\), confirming that the EnKF continuation is selected from the interior of the offline assimilation window rather than from \(T\).
The observation and correction contributions to the local loss, together with their bell-averaged value used for restart selection, are shown in \Cref{f:l96_restart_loss_representative}.
Across all twenty experiments, the mean loss-selected restart-state RMSE is \(0.180\), compared with a mean RMSE of \(0.438\) at \(T\), corresponding to a \(59\%\) reduction.
The restart-state RMSE is lower than the RMSE at \(T\) in sixteen of the twenty experiments.

\begin{figure}[!ht]
\centering
\includegraphics[width=0.8\textwidth]{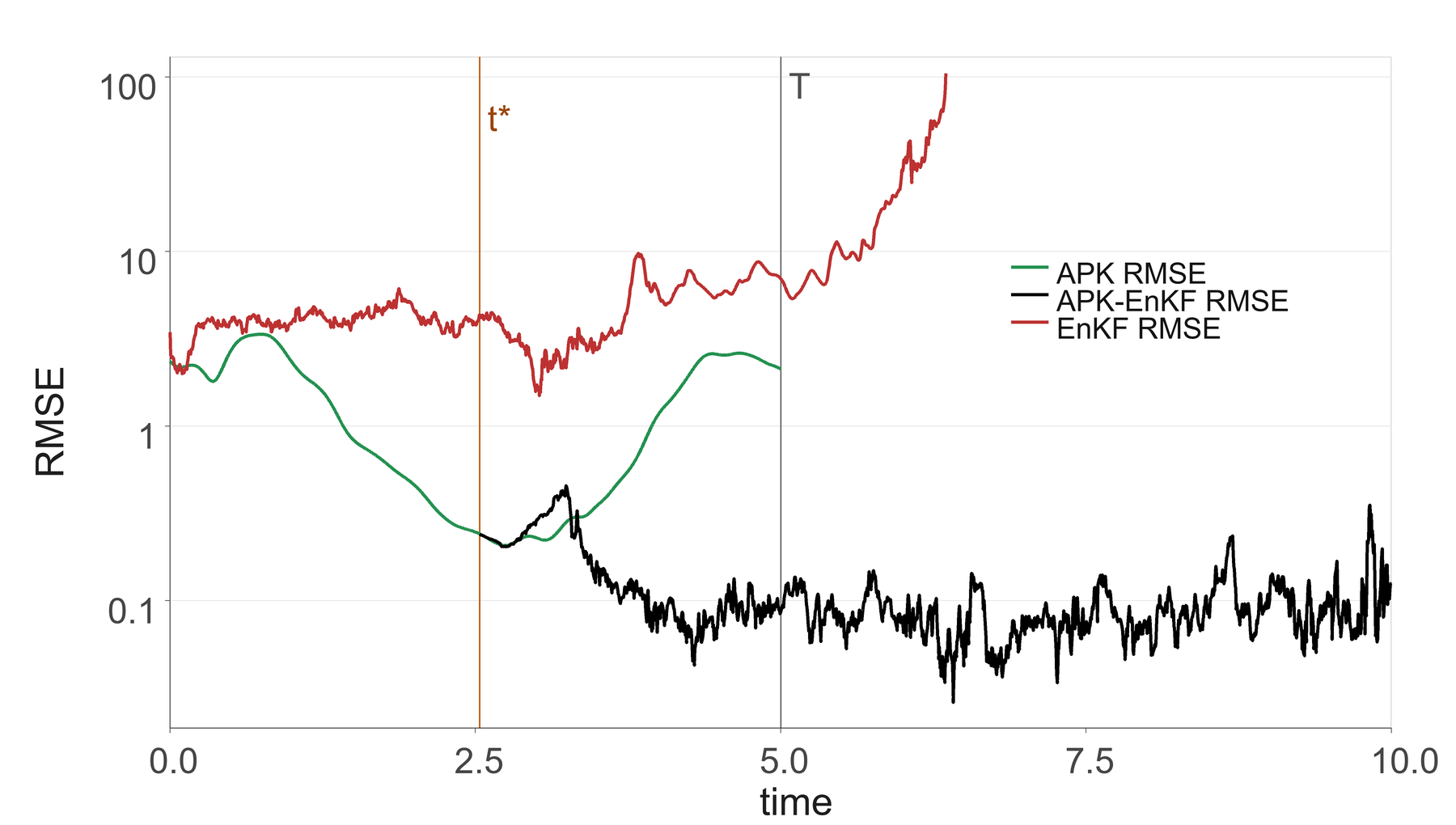}
\caption{
Representative all-coordinate RMSE over physical time.
The green curve is the APK reconstruction on \(0\leq t\leq T\), the black curve is the EnKF initialized at the loss-selected time \(t_*\), and the red curve is the ordinary EnKF from the rough prior.
The gray line marks \(T=5\); the orange line marks \(t_*\approx2.54\).
The ordinary-EnKF curve ends when its state becomes nonfinite.
}
\label{f:l96_filter_representative}
\end{figure}

\begin{figure}[!ht]
\centering
\includegraphics[width=\textwidth]{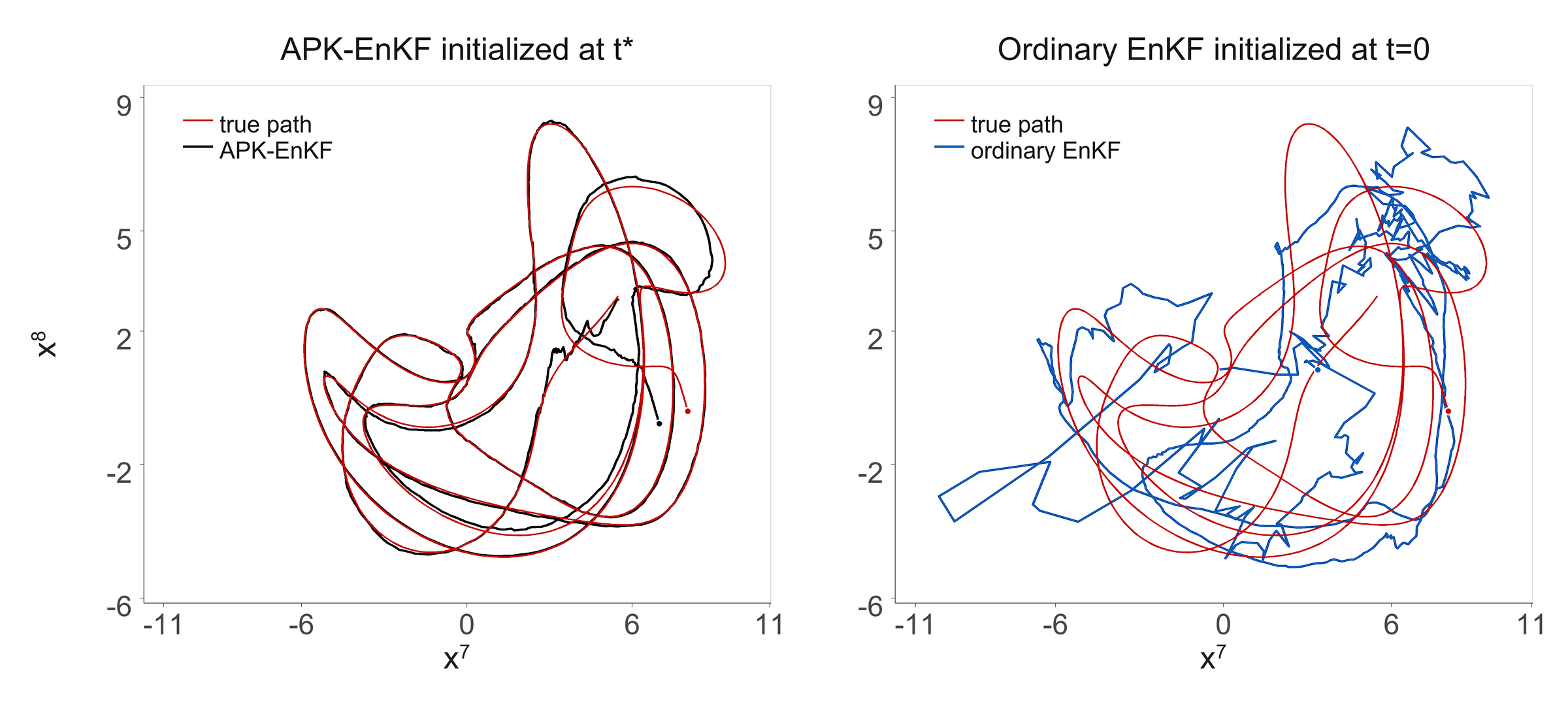}
\caption{
Representative trajectories in the two hidden coordinates \((x^7,x^8)\).
Both panels show \(t_*\leq t\leq2T\), and the dots mark the states at \(t_*\).
The left panel compares the true path with APK--EnKF initialized at \(t_*\), while the right panel compares the true path with ordinary EnKF initialized at \(t=0\), with only the ordinary-EnKF trajectory after \(t_*\) displayed.
The APK-initialized filter remains near the selected branch, whereas the ordinary EnKF follows a different high-error trajectory.
}
\label{f:l96_hidden_trajectory_representative}
\end{figure}

Together, the aggregate and representative results separate global branch selection from local online tracking.
With the same locally stable EnKF tuning, the APK restart leads to accurate continuation in every experiment, while the rough prior either remains far from the truth or loses numerical stability.

\FloatBarrier
\subsection{Separate nonlinear-observation setting}
\label{s:numeric_squared}

The separate squared-observation setting tests whether the offline advantage persists under a nonlinear observation map.
We observe one squared coordinate and then hide the next three:
\[
\phi(x)=(x_0^2,x_4^2,\ldots,x_{36}^2).
\]
The system dimension, noise standard deviation \(0.3\), optimization budgets, and parameter schedules are unchanged, while the correction penalty is recalibrated to \(C\approx0.540\).
Across twenty independent experiments, APK has a mean interior path RMSE of \(1.91\pm0.692\), compared with \(3.73\pm0.543\) for population weak-\(\mathrm{4D\text{-}Var}_x\).
Their loss-selected restart-state RMSEs are \(0.864\pm0.525\) and \(3.30\pm0.941\), respectively.

The nonlinear observation map makes continuation less uniform.
APK--EnKF has an online RMSE below \(0.5\) in fifteen of twenty experiments and a median RMSE of \(0.0134\), but its mean is \(1.68\pm2.96\) because five continuations select poor branches.
Ordinary EnKF becomes nonfinite in one experiment and has a mean RMSE of \(6.64\pm0.391\) among its nineteen complete experiments.
Thus the squared-observation setting preserves the offline advantage, but in five experiments the selected run appears to settle in a local basin whose loss-selected state does not support accurate continuation.
This motivates adapting \(\alpha\) and \(\sigma\) in time and space to renew exploration where optimization stagnates.

\FloatBarrier
\section{Conclusion and future work}
\label{s:conclusion}

This paper develops a full-window framework for offline branch discovery and optional online EnKF continuation.
The branch search combines stable APK gradients, an optimized random initial law, prescribed path-kernel and diffusion schedules, and population mixing across independent runs.
In the direct long-window implementations tested here, strong-constraint \(\mathrm{4D\text{-}Var}\) and weak-\(\mathrm{4D\text{-}Var}_q\) have nonfinite gradients at the first optimization step.
Against the remaining finite variational comparator, population weak-\(\mathrm{4D\text{-}Var}_x\), APK finds more accurate paths and restart states.
A computable local loss selects an interior APK state from which the EnKF may continue online; the resulting APK--EnKF remains accurate when ordinary EnKF is inaccurate or becomes nonfinite.

Adaptive exploration across physical time and state coordinates is a natural next step.
Figure~\ref{f:l96_restart_loss_representative} suggests a truth-free signal for this adaptation: large local loss tends to coincide with large diagnostic RMSE.
Although the agreement is not exact, the local spatiotemporal loss could therefore help identify, without access to the true state, which physical times and state coordinates should use larger \(\alpha\) and \(\sigma\).
Locations with small loss could instead use smaller values and emphasize local exploitation.
The APK theory permits adapted path-kernel schedules and nonconstant diffusion, so this extension remains within the same theoretical framework.
By concentrating exploration where the current reconstruction is least satisfactory, it may improve branch discovery without adding unnecessary randomness elsewhere.

\section*{Data availability statement}

The code and generated numerical data supporting this study are available from the author upon reasonable request.

\bibliographystyle{abbrv}
{\footnotesize\bibliography{library}}

\begin{thebibliography}{10}

\bibitem{Ahrens1999}
B.~Ahrens.
\newblock Variational data assimilation for a lorenz model using a non-standard
  genetic algorithm.
\newblock {\em Meteorology and Atmospheric Physics}, 70(3--4):227--238, 1999.

\bibitem{filter_ABN}
M.~Asch, M.~Bocquet, and M.~Nodet.
\newblock {\em Data Assimilation}.
\newblock Society for Industrial and Applied Mathematics, 12 2016.

\bibitem{AurouxBlum2008}
D.~Auroux and J.~Blum.
\newblock A nudging-based data assimilation method: The back and forth nudging
  ({BFN}) algorithm.
\newblock {\em Nonlinear Processes in Geophysics}, 15(2):305--319, 2008.

\bibitem{AzouaniOlsonTiti2014}
A.~Azouani, E.~Olson, and E.~S. Titi.
\newblock Continuous data assimilation using general interpolant observables.
\newblock {\em Journal of Nonlinear Science}, 24(2):277--304, 2014.

\bibitem{BocquetSakov2014IEnKS}
M.~Bocquet and P.~Sakov.
\newblock An iterative ensemble {Kalman} smoother.
\newblock {\em Quarterly Journal of the Royal Meteorological Society},
  140:1521--1535, 2014.

\bibitem{Brocker2017}
J.~Bröcker.
\newblock Existence and uniqueness for four-dimensional variational data
  assimilation in discrete time.
\newblock {\em SIAM Journal on Applied Dynamical Systems}, 16:361--374, 1 2017.

\bibitem{Caflisch2021}
R.~Caflisch, D.~Silantyev, and Y.~Yang.
\newblock Adjoint dsmc for nonlinear boltzmann equation constrained
  optimization.
\newblock {\em Journal of Computational Physics}, 439:110404, 8 2021.

\bibitem{CM44}
R.~H. Cameron and W.~T. Martin.
\newblock Transformations of weiner integrals under translations.
\newblock {\em The Annals of Mathematics}, 45:386, 4 1944.

\bibitem{Cox2015}
G.~Cox.
\newblock (non)uniqueness of critical points in variational data assimilation.
\newblock {\em Physica D: Nonlinear Phenomena}, 300:34--40, 4 2015.

\bibitem{4dvar_DT}
F.-X.~L. DIMET and O.~TALAGRAND.
\newblock Variational algorithms for analysis and assimilation of
  meteorological observations: theoretical aspects.
\newblock {\em Tellus A}, 38A:97--110, 3 1986.

\bibitem{eyink2004ruelle}
G.~L. Eyink, T.~W.~N. Haine, and D.~J. Lea.
\newblock Ruelle's linear response formula, ensemble adjoint schemes and lévy
  flights.
\newblock {\em Nonlinearity}, 17:1867--1889, 9 2004.

\bibitem{Glynn1990}
P.~W. Glynn.
\newblock Likelihood ratio gradient estimation for stochastic systems.
\newblock {\em Communications of the ACM}, 33:75--84, 10 1990.

\bibitem{DA_Kalnay}
E.~Kalnay, S.~Mote, and C.~Da.
\newblock {\em Earth System Modeling, Data Assimilation and Predictability}.
\newblock Cambridge University Press, 10 2024.

\bibitem{filter_LDZ}
Z.~Li, B.~Dong, and P.~Zhang.
\newblock Latent assimilation with implicit neural representations for unknown
  dynamics.
\newblock {\em Journal of Computational Physics}, 506:112953, 6 2024.

\bibitem{lucarini_linear_response_climate}
V.~Lucarini, F.~Ragone, and F.~Lunkeit.
\newblock Predicting climate change using response theory: Global averages and
  spatial patterns.
\newblock {\em Journal of Statistical Physics}, 166:1036--1064, 2017.

\bibitem{MalliavinBook}
P.~Malliavin.
\newblock {\em Stochastic Analysis}, volume 313.
\newblock Springer Berlin Heidelberg, 1997.

\bibitem{filter_NMX}
A.~Narayan, Y.~Marzouk, and D.~Xiu.
\newblock Sequential data assimilation with multiple models.
\newblock {\em Journal of Computational Physics}, 231:6401--6418, 8 2012.

\bibitem{apk}
A.~Ni.
\newblock Adjoint path-kernel method for backpropagation and data assimilation
  in unstable diffusions.
\newblock {\em arXiv:2507.21497}, 7 2025.

\bibitem{dud}
A.~Ni.
\newblock Differentiating unstable diffusion.
\newblock {\em arXiv:2503.00718}, 3 2025.

\bibitem{PoterjoyZhang2015}
J.~Poterjoy and F.~Zhang.
\newblock Systematic comparison of four-dimensional data assimilation methods
  with and without the tangent linear model using hybrid background error
  covariance: {E4DVar} versus {4DEnVar}.
\newblock {\em Monthly Weather Review}, 143:1601--1621, 2015.

\bibitem{Reiman1989}
M.~I. Reiman and A.~Weiss.
\newblock Sensitivity analysis for simulations via likelihood ratios.
\newblock {\em Operations Research}, 37:830--844, 10 1989.

\bibitem{Rubinstein1989}
R.~Y. Rubinstein.
\newblock Sensitivity analysis and performance extrapolation for computer
  simulation models.
\newblock {\em Operations Research}, 37:72--81, 2 1989.

\bibitem{DA_AST}
D.~Sanz-Alonso, A.~Stuart, and A.~Taeb.
\newblock {\em Inverse Problems and Data Assimilation}.
\newblock Cambridge University Press, 7 2023.

\bibitem{Sanz-Alonso2025}
D.~Sanz-Alonso and N.~Waniorek.
\newblock Long-time accuracy of ensemble kalman filters for chaotic dynamical
  systems and machine-learned dynamical systems.
\newblock {\em SIAM Journal on Applied Dynamical Systems}, 24:2246--2286, 9
  2025.

\bibitem{4dvar_TC}
O.~Talagrand and P.~Courtier.
\newblock Variational assimilation of meteorological observations with the
  adjoint vorticity equation. i: Theory.
\newblock {\em Quarterly Journal of the Royal Meteorological Society},
  113:1311--1328, 10 1987.

\bibitem{Tremolet2006}
Y.~Tr{\'e}molet.
\newblock Accounting for an imperfect model in 4d-var.
\newblock {\em Quarterly Journal of the Royal Meteorological Society},
  132:2483--2504, 2006.

\bibitem{XuAnitescu2016}
W.~Xu and M.~Anitescu.
\newblock A limited-memory multiple shooting method for weakly constrained
  variational data assimilation.
\newblock {\em SIAM Journal on Numerical Analysis}, 54:3300--3331, 2016.

\bibitem{filter_Yau}
S.~S.-T. Yau, X.~Chen, X.~Jiao, J.~Kang, Z.~Sun, and Y.~Tao.
\newblock {\em Principles of nonlinear filtering theory}.
\newblock SPRINGER INTERNATIONAL PU, 2024.

\bibitem{ZhangZhang2012E4DVar}
M.~Zhang and F.~Zhang.
\newblock {E4DVar}: Coupling an ensemble kalman filter with four-dimensional
  variational data assimilation in a limited-area weather prediction model.
\newblock {\em Monthly Weather Review}, 140:587--600, 2012.

\bibitem{ZhaoWangLiu2011}
J.~Zhao, B.~Wang, and J.~Liu.
\newblock Impact of analysis-time tuning on the performance of the {DRP-4DVar}
  approach.
\newblock {\em Advances in Atmospheric Sciences}, 28:207--216, 2011.

\end{thebibliography}

\end{document}